\RequirePackage{fix-cm}
\documentclass[smallextended]{svjour3}       
\smartqed  
\usepackage{graphicx}
\usepackage{a4}
\usepackage{graphicx}
\usepackage{parskip}

\usepackage{amsfonts}
\usepackage{natbib}
\usepackage{algorithm}
\usepackage{algorithmic}

\usepackage{latexsym}

\newcommand{\R}{\mathbb{R}}
\newcommand{\Nat}{\mathbb{N}}

\newcommand{\PP} {{  \rm I\hskip-0.22em P}}

\newcommand{\EE} {{\rm I\hskip-0.48em E}}

\usepackage{natbib}

\begin{document}

\title{The Bernstein-Orlicz norm and deviation inequalities\thanks{
Research supported by SNF 20PA21-120050.}}


\author{Sara van de Geer         \and
        Johannes Lederer 
}


\institute{S. van de Geer \at
              Seminar for Statistics, ETH Z\"urich\\ R\"amistrasse 101, 8092 Z\"urich, Switzerland\\
              Tel.:  00-41-44-6322252\\
              Fax:  00-41-44-6321228\\
              \email{geer@stat.math.ethz.ch}           
           \and
           J. Lederer \at
           Seminar for Statistics, ETH Z\"urich\\ R\"amistrasse 101, 8092 Z\"urich, Switzerland\\
}

\date{Received: date / Accepted: date}

\maketitle

\begin{abstract}
We introduce two new concepts designed for the study of empirical
processes. First, 
we introduce a new Orlicz norm which we call the Bernstein-Orlicz norm.
This new norm interpolates sub-Gaussian and
sub-exponential tail behavior. In particular, we show how this norm can be 
used to simplify the derivation of  deviation inequalities for suprema of collections of
random variables.
Secondly, we introduce chaining and generic chaining along a tree.
These simplify the well-known concepts of chaining and generic chaining.
The supremum of  the empirical process is then
studied as a special case. We show that chaining along a tree can be
done using entropy with bracketing. Finally, we establish a deviation inequality
for the empirical process for the unbounded case.
\end{abstract}

\section{Introduction}\label{introduction.section}
We introduce a new Orlicz norm which we name the Bernstein-Orlicz norm.
It interpolates sub-Gaussian and
sub-exponential tail behavior. With this new norm, we apply the
usual techniques based on Orlicz norms. In particular, we derive
deviation
inequalities for suprema in a fairly
simple and straightforward way.
The Bernstein-Orlicz norm captures Bernstein's
probability inequalities, and its use puts further derivations in a unifying framework,
shared for example by techniques for the  sub-Gaussian case,
such as those for empirical processes based on symmetrization and Hoeffding's inequality.

We furthermore introduce chaining and generic chaining along a tree, which is
we believe conceptually simpler than the usual chaining and generic chaining. 
We invoke it for the presentation of maximal inequalities for general random variables with
finite Bernstein-Orlicz norm. The supremum of  the empirical process is then
studied as a special case, and we show that chaining along a tree can be
done using entropy with bracketing. We establish a deviation inequality
for the empirical process indexed by a class of functions
${\cal G}$, in terms of the new Bernstein-Orlicz norm.
The class ${\cal G}$ is assumed to satisfy a uniform Bernstein
condition, but need not be uniformly bounded in supremum norm.

The paper is organized as follows. In Section
 \ref{definition.section}, we introduce
the Bernstein-Orlicz norm and discuss the relation with
Bernstein's inequality. We then present some bounds for
maxima of finitely many random variables (Section \ref{finitely-many.section}) or suprema 
over a countable set of random variables (Section \ref{Suprema.section}). Section \ref{Suprema.section} also contains the concept of
(generic) chaining along a tree. The proofs of the results in Sections \ref{definition.section}, 
\ref{finitely-many.section} and \ref{Suprema.section}
are elementary and given immediately following their statement.
Section \ref{adaptivetruncationsection} contains the application to the
empirical process. The proofs here are more technical, and given separately
in Sections \ref{proofs.section} and \ref{lastsection}.

\section{The Bernstein-Orlicz norm} \label{definition.section}

Consider a random variable $Z \in \R$ with distribution $\PP$.
We first recall the general Orlicz norm (see e.g.\ \cite {krasnosel1961convex}).

\begin{definition} Let $\Psi: [0, \infty) \rightarrow
[0 , \infty ) $ be an increasing and convex function with $\Psi(0)=0$. The 
$\Psi$-Orlicz norm of $Z$ is 
$$\| Z \|_{\Psi} := \inf \biggl  \{ c >0 : \ \EE \Psi \biggl ({ |Z | \over c }  \biggr ) \le 1 \biggr \} . $$
\end{definition}

A special case is the $L_m (\PP)$-norm ($m \ge 1$)
which corresponds to $\Psi (z) = z^m $. Other important special cases 
are $\Psi (z) = \exp [ z^2 ] - 1$ for sub-Gaussian random variables and
$\Psi(z) = \exp (z) -1 $ for sub-exponential random variables.
We propose functions $\Psi$ that combine sub-Gaussian intermediate tails 
and sub-exponential far tails.

For each $L>0$ we define 
\begin{equation}\label{PsiL.equation}
\Psi_L (z):= \exp \biggl [ { \sqrt {1+ 2L z } - 1 \over   L  } \biggr ]^2 - 1 , \ z \ge 0 .  
\end{equation}
It is easy to see that $\Psi_L$ is increasing and convex, and that
$\Psi_L (0)=0$.

\begin{definition} Let $L>0$ be given. The ($L$-)Bernstein-Orlicz norm is the $\Psi$-Orlicz norm with $\Psi= \Psi_L$ given in
(\ref{PsiL.equation}). 
\end{definition}

Indeed, the Bernstein-Orlicz norm combines sub-Gaussian and
sub-exponential behavior:
$$\Psi_L (z) \approx \cases {\exp [ z^2] -1 & for $Lz$ small \cr
\exp[2z/L ] -1 &  for $Lz$ large \cr  } . $$

Note that the constant $L$ governs the range of the sub-Gaussian behavior.
It is a dimensionless constant, i.e., it does not depend on the scale of measurement.

The inverse of $\Psi_L$ is
$$\Psi_L^{-1} (t)= \sqrt { \log (1+t)} + {L \over 2} \log (1+t) , \ t \ge 0 . $$
With this and with Chebyshev's inequality, one now directly derives a probability
inequality for $Z$. 

\begin{lemma} \label{Psi-Prob.lemma} Let
$\tau:= \| Z \|_{\Psi_L} $.  We have for all $t >0$, 
$$\PP \biggl ( |Z|  >\tau \biggl [  \sqrt {t} + {Lt \over 2} \biggr ]  \biggr ) \le 2 \exp [-t] . $$
\end{lemma}

{\bf Proof of Lemma \ref{Psi-Prob.lemma}.} By Chebyshev's inequality, for all $c>\| Z \|_{\Psi_L} $, 
$$\PP \biggl ( | Z | / c  \ge  \sqrt {t} + {Lt \over 2} \biggr ) = 
\PP \biggl ( |Z| / c \ge \Psi_L^{-1} ( {\rm e}^t -1) \biggr ) $$
$$ = \PP \biggl ( \Psi_L (| Z|/c)  \ge {\rm e}^t -1  \biggr ) \le 
\biggl ( \EE \Psi_L ( |Z|/ c ) +1 \biggr ) {\rm e}^{-t} .$$
Thus,
$$\PP \biggl ( | Z | / \tau  >  \sqrt {t} + {Lt \over 2} \biggr ) =
\lim_{c \downarrow \tau}
\PP \biggl ( | Z | / c  >  \sqrt {t} + {Lt \over 2} \biggr ) $$
$$ \le \lim_{c \downarrow \tau} 
\biggl ( \EE \Psi_L ( |Z|/ c ) +1 \biggr ) {\rm e}^{-t} \le 2 {\rm e}^{-t} . $$

\hfill $\sqcup \mkern -12mu \sqcap$

The next lemma says that a converse result holds as well, 
that is, from the probability inequality of Lemma 
\ref{Psi-Prob.lemma} one can derive a bound for the Bernstein-Orlicz
norm, with constants $L$ and $\tau$ multiplied by $\sqrt 3$\footnote{The constant
can possibly be improved.}.

 \begin{lemma} \label{Prob-Psi.lemma} Suppose that for for some constants
 $\tau$ and $L$, and for all $t >0$,
  $$\PP \biggl ( |Z| \ge \tau \biggl [ \sqrt {t} + {Lt  \over 2}  \biggr ] \biggr ) \le
  2 \exp[-t] . $$
  Then
  $\| Z \|_{\Psi_{\sqrt 3 L }} \le \sqrt 3 \tau $. 
  \end{lemma}
  
  {\bf Proof of Lemma \ref{Prob-Psi.lemma}.} We have
  $$\EE \Psi_{\sqrt 3L} \biggl (|Z|/(\sqrt {3} \tau)\biggr )=
  \int_0^{\infty} \PP \biggl ( |Z| \ge \sqrt 3 \tau \Psi_{\sqrt 3 L}^{-1} (t) \biggr )dt $$
  $$ = \int_0^{\infty} \PP \biggl ( |Z| \ge  \sqrt 3 \tau 
  \biggl [ \sqrt {\log (1+t) } + {\sqrt 3 L \over 2 } \log (1+t ) \biggr ] \biggr )dt  $$ 
  $$ = \int_0^{\infty} \PP \biggl ( |Z| \ge   \tau 
  \biggl [ \sqrt {\log (1+t)^3 } + { L \over 2 } \log (1+t )^3 \biggr ] \biggr )dt  
      \le 2 \int_{0}^{\infty} {1 \over (1+t)^3} dt =  1.
  $$
  \hfill $\sqcup \mkern -12mu \sqcap$

We recall Bernstein's inequality, see \cite{Bennet:62}. 

\begin{theorem} \label{Prob-Bernstein.theorem}
Let $X_1 , \ldots , X_n $ be independent random variables
with values in $\R$ and with mean zero. Suppose that for some constants $\sigma$ and $K$,
one has
$$
{1 \over n} \sum_{i=1}^n \EE | X_i|^m \le {m! \over 2} K^{m-2} \sigma^2 , \ m =1 , 2 , \ldots . $$
Then for all $t>0$,
$$ \PP \biggl ( {1 \over \sqrt n } \biggl | \sum_{i=1}^n X_i \biggr | \ge \sigma \sqrt {2 t } + { Kt \over\sqrt  n} \biggr ) \le 
2 \exp [-t] . $$
\end{theorem}

The following corollary shows that $\| \cdot \|_{\Psi_L}$ indeed captures the nature
of Bernstein's inequality.

\begin{corollary} \label{Bernstein-Psi.corollary}
Let $X_1 , \ldots , X_n $ be independent random variables
satisfying the conditions of Theorem \ref{Prob-Bernstein.theorem}.
Then by this theorem and Lemma \ref{Prob-Psi.lemma}, 
for $L := \sqrt 6 K /( \sqrt n \sigma)  $, we have
$$\biggl \| {1 \over \sqrt n } \sum_{i=1}^n X_i \biggr \|_{\Psi_L} \le  \sqrt 6 \sigma . $$

\end{corollary}

\section{The Bernstein-Orlicz norm for the maximum of finitely many variables}
\label{finitely-many.section}

Using Orlicz norms, the argument for obtaining
a bound for the expectation of maxima is standard.
We refer to \cite{vanderVaart:96} for a general approach. 
We consider the special case of the Bernstein-Orlicz norm.

\begin{lemma}\label{Expect-Finite.lemma}
Let $\tau$ and $L$ be constants, and let $Z_1 , \ldots , Z_p$ be random variables satisfying
$$\max_{1 \le j \le p } \| Z_j \|_{\Psi_L} \le \tau . $$
Then
$$\EE  \max_{1 \le j \le p} |Z_j | \le \tau\biggl [ \sqrt {\log(1+ p)} + {L\over 2}  \log (1+p )
\biggr ]  . $$
\end{lemma}

{\bf Proof of Lemma \ref{Expect-Finite.lemma}
.} Let $c > \tau$. Then by Jensen's inequality
$$\EE  \max_{1 \le j \le p} |Z_j |\le c \Psi_L^{-1} \left  ( \EE \Psi_L\biggl (
\max_{1 \le j \le p } |Z_j | / c\biggr  ) \right )
 = c \Psi_L^{-1} \left  ( \EE \max_{1 \le j \le p } \Psi_L\biggl (
 |Z_j | / c\biggr  ) \right )  $$ $$ \le c \Psi_L^{-1} \left  ( \sum_{j=1}^p 
 \EE \Psi_L\biggl (
 |Z_j | / c\biggr  ) \right ) \le c \Psi_L^{-1} \left ( p \max_{1 \le j \le p }
  \EE \Psi_L\biggl (
 |Z_j | / c\biggr  ) \right ). 
 $$
 Therefore,
 $$\EE  \max_{1 \le j \le p} |Z_j | \le \lim_{c \downarrow  \tau} 
  c \Psi_L^{-1} \left ( p \max_{1 \le j \le p }
  \EE \Psi_L\biggl (
 |Z_j | / c\biggr  ) \right ) \le \tau \Psi_L^{-1} ( p ) $$ $$=
  \tau \biggl [ \sqrt {\log(1+ p)} + {L\over 2}  \log (1+p ) \biggr ] . $$
  \hfill $\sqcup \mkern -12mu \sqcap$
  
  As a special case, one may consider the random variables
  $$Z_j := {1 \over \sqrt n} \sum_{i=1}^n  g_j(X_i), \ j=1 , \ldots , p , $$
  where $X_1 , \ldots , X_n$ are independent random variables with values
  in some space ${\cal X}$, and where $g_1 , \ldots , g_p$ are real-valued functions on
  ${\cal X}$.  If the $g_j(X_i)$ are centered for all $i$ and
  $j$, and if one assumes the Bernstein condition
  $${1 \over n } \sum_{i=1}^n \EE | g_j (X_i) |^m \le { m! \over 2 } K^{m-2} \sigma^2 , m=2,3, \ldots , \ 
  j=1 , \ldots, p , $$
  then one can apply Lemma \ref{Expect-Finite.lemma}, with
  $\tau:= \sqrt 6 \sigma$ and $L= \sqrt 6 K/ (\sqrt n \sigma)$, giving the inequality
  \begin{equation} \label{S=0}
  \EE \max_{1 \le j \le p} \biggl | {1 \over \sqrt n}  \sum_{i=1}^n g_j (X_i)  \biggr | 
  \le \sigma \sqrt {6 \log (1+ p)} + { 3K \over \sqrt n} \log (1+ p ) . 
  \end{equation}
   This follows from
  Corollary \ref{Bernstein-Psi.corollary}. The constants can however be improved when
  using direct arguments (see e.g.\ Lemma 14.12 \cite{BvdG2011}). 
  
  We now present a deviation inequality in probability for the maximum of
  finitely many variables.
  \begin{lemma}\label{deviation-Prob.lemma}
  Let Let $Z_1 , \ldots , Z_p$ be random variables satisfying for some $L$ and
$\tau$
$$\max_{1 \le j \le p } \| Z_j \|_{\Psi_L} \le \tau . $$
Then for all $t > 0$
$$\PP \left ( \max_{1 \le j \le p } | Z_j | \ge 
  \tau \biggl [ \sqrt {\log(1+ p)} + {L\over 2}  \log (1+p ) + \sqrt t +{ Lt \over 2}   \biggr ] 
  \right ) 
  \le 2 \exp[-t] . $$
  \end{lemma}
  
  {\bf Proof of Lemma \ref{deviation-Prob.lemma}.} We first use that for any $a>0$ and $t>0$, one has $\sqrt {a} + \sqrt {t} >
  \sqrt {a+t}$, so that
  $$\PP \left ( \max_{1 \le j \le p } | Z_j | \ge  
  \tau \biggl [ \sqrt {\log(1+ p)} + {L\over 2}  \log (1+p ) + \sqrt t +{ Lt \over 2}   \biggr ] 
  \right ) $$
$$ \le   \PP \left ( \max_{1 \le j \le p } | Z_j | > 
  \tau \biggl [ \sqrt {t+\log(1+ p)} + {L\over 2} ( t+ \log (1+p ) ) \biggr ] 
  \right ) .$$
  Next, we apply the union bound and Lemma \ref{Psi-Prob.lemma}:
 $$  \PP \left ( \max_{1 \le j \le p } | Z_j | > 
  \tau \biggl [ \sqrt {t+\log(1+ p)} + {L\over 2} ( t+ \log (1+p ) ) \biggr ]
  \right ) $$ $$ \le \sum_{j=1}^p 
  \PP \left ( | Z_j | > 
  \tau \biggl [ \sqrt {t+\log(1+ p)} + {L\over 2} ( t+ \log (1+p ) ) \biggr ] 
  \right )$$
  $$ \le 2 p \exp\biggl [ -(t+ \log (1+p) ) \biggr ] = {2p \over 1+p} \exp[-t] \le 2 \exp[-t] . $$
  \hfill $\sqcup \mkern -12mu \sqcap$
  
  Using Lemma \ref{Prob-Psi.lemma}, this is easily converted into a the
  following deviation
  inequality for the Bernstein-Orlicz norm.
  We use the notation
  $$x_+ := x {\rm l } \{ x > 0 \} . $$

  \begin{lemma}\label{deviation-Psi.lemma}
   Let Let $Z_1 , \ldots , Z_p$ be random variables satisfying for some $L$ and
$\tau$
$$\max_{1 \le j \le p } \| Z_j \|_{\Psi_L} \le \tau . $$
Then 
$$\biggl \| \biggl ( \max_{1 \le j \le p } 
|Z_j| - \tau \biggl [  \sqrt {\log (1+p)}  + {L \over 2} \log (1+p)  \biggr ]  \biggr )_+ \biggr \| _{\Psi_{\sqrt 3 L}} \le \sqrt 3 \tau. $$
\end{lemma}
  
  {\bf Proof of Lemma \ref{deviation-Psi.lemma}.} Let 
 $$Z := \biggl  (\max_{1 \le j \le p } |Z_j| -  \tau \biggl [    \sqrt {\log(1+ p)} + {L\over 2}  \log (1+p ) \biggr ]  \biggr )_+ . $$
 By Lemma \ref{deviation-Prob.lemma}, we have for all $t >0$
 $$ \PP \biggl ( Z \ge  \tau \biggl [ \sqrt t + {Lt \over 2} \biggr ] \biggr )$$ $$ =
 \PP \left ( \max_{1 \le j \le p } | Z_j | \ge 
  \tau \biggl [ \sqrt {\log(1+ p)} + {L\over 2}  \log (1+p ) + \sqrt t +{ L t\over 2}   \biggr ]
  \right ) 
  \le 2 \exp[-t] . $$
 Application of Lemma \ref{Prob-Psi.lemma} finishes the proof.
 
 \hfill $\sqcup \mkern -12mu \sqcap$.

   \section{Chaining along a tree}\label{Suprema.section}

  A common technique for bounding suprema of stochastic processes is 
  chaining as developed by Kolmogorov, leading to versions of 
  Dudley's entropy bound (\cite{dudley1967sizes}). See e.g. \cite {vanderVaart:96} or 
  \cite{vandeGeer:00} and the references therein.
  We however propose another method which we call chaining along a tree.  This
  method  is 
  conceptually simpler than the usual chaining and,
as far as we know, does not introduce unnecessary restrictions. 
An example will be detailed in Section \ref{adaptivetruncationsection}
  for the case of entropy with bracketing.
The generic chaining technique of \cite{talagrand2005generic} is a refinement which we
shall also consider in Definition \ref{generic-tree.definition} and
Theorem \ref{generic}. 

Let $S \in \Nat_0$ be fixed.

\begin{definition}  A finite 
tree\footnote{Actually, ${\cal T}$ is rather a forest consisting of
$|G_0|$ trees} ${\cal T} $ is a collection $\{ G_s \}_{s=0}^S $
of disjoint subsets of $\{ 1 , \ldots , N\} $  such that
$\cup_{s=0}^S G_s = \{ 1 , \ldots , N \}$, together with a function
$${\rm parent} : \{ 1 , \ldots , N\} \rightarrow \{ 1 , \ldots , N\} ,$$ such that
${\rm parent } (j) \in G_{s-1}$ for $j \in G_s $, $s \in 1 , \ldots , S $. We call
an element of 
$\{ 1 , \ldots , N\}$ a node, and $G_s$ a generation, $s=0 , \ldots , S$. 
A branch of the tree with end node $j_S \in G_S$ is the sequence
$\{ j_0 , \ldots , j_S\} $ with $j_{s-1} = {\rm parent}({j_s})$, $s=1, \ldots , S$.
\end{definition}

\begin{definition} Let a collection of real-valued random
variables ${\cal W }:= \{ W_j \}_{j=1}^N$ be given.
A finite labeled tree $({\cal T}, {\cal W})$ is a finite tree
with on each node $j$ a label $W_j $.
\end{definition}

 Let $\Theta$ be some countable set and let
  $Z_{\theta} \in \R$ be a random variable defined for
  each $\theta \in \Theta$. 
  We consider supremum of the process $\{ |Z_{\theta} | : \ \theta \in \Theta \} $. 

\begin{definition} \label{tree.definition}
Let $\delta >0$ and $\tau>0$ be constants and let ${\cal L} := \{ L_s \}_{s=0}^S$ be a sequence
of positive numbers.
A $(\delta,\tau,  {\cal L} )$ finite tree chain for $\{ Z_{\theta} \}$ is a finite labeled
tree $({\cal T} ,{\cal W} )$ such that for all $s=0 , \ldots , S$, 
$$ \| W_{j}  \|_{\Psi_{L_s}} \le \tau 2^{-s} , \ \forall \ j \in G_s, $$
and such that one can apply chaining of $\{ Z_{\theta } \}$ along the
tree $({\cal T} , {\cal W})$, with approximation error
$\delta$. That is, 
for each $\theta \in \Theta$ there is an end  node $j_S \in
G_S$ such that the branch $\{ j_0 , \ldots , j_S \}$ satisfies
$$| Z_{\theta}| \le   \sum_{s=0}^S | W_{j_s} | + \delta . $$
\end{definition}

In the above definition, the approximation error $\delta$ will generally
depend on the depth $S$ of the tree. We assume that at a
fine enough level, the approximation error is small.
The usual chaining technique does not assume a tree
structure, but indeed often needs only a finite
number of steps. A tree structure follows if the members at the finest
level are taken as end nodes. With a finite number of steps,
the sum given in (\ref{gamma.equation}) is finite. 
This avoids requiring  convergence of an infinite sum. 

We have presented the definition of a finite tree chain for the
Bernstein-Orlicz norm $\| \cdot \|_{\Psi_L}$.  However, the concept is
not particularly tied up with this norm, e.g., for sub-Gaussian cases
one may choose to replace the Bernstein-Orlicz norm by
the $L_2 (\PP)$ norm (corresponding to  case where the constants in ${\cal L}$
all vanish). 

Let us now turn to the results.

\begin{theorem} \label{Tree.theorem}
Let $({\cal T},{\cal W} )$ be an $(\delta ,\tau, {\cal  L} )$ finite tree chain for
$\{ Z_{\theta } \}$. Define
\begin{equation}\label{gamma.equation}
\gamma:= \tau \sum_{s=0}^S 2^{-s}
\biggl [ \sqrt {\log (1+ |G_s| ) } + {L_s \over 2} \log (1+ |G_s| )  \biggr ] . 
\end{equation}
It holds that
\begin{equation}\label{meanbound.equation}
\EE \left ( \sup_{\theta \in \Theta} | Z_{\theta} | \right ) \le 
  \gamma + \delta. 
  \end{equation}
\end{theorem}

\begin{remark}\label{minimize-tree.remark}
One
may minimize the right hand side of (\ref{meanbound.equation}) over all 
finite trees. 
\end{remark}

{\bf Proof of Theorem \ref{Tree.theorem}.} We have
$$\EE \sup_{\theta \in \Theta} |Z_{\theta} | \le
  \sum_{s=0}^S \EE \max_{j \in G_s } | W_j |+  \delta .$$
Application of Lemma \ref{Expect-Finite.lemma} gives that for each
$s \in \{ 0 , \ldots , S \} $
$$\EE \max_{j \in G_s } |W_j| \le \tau 2^{-s} \biggl [ \sqrt {\log (1+ |G_s| ) } +
{ L_s  \over 2}  \log (1+ | G_s | ) \biggr ] . $$
\hfill $\sqcup \mkern -12mu \sqcap$

With generic chaining, the condition on the Bernstein-Orlicz norm
of the labels is dropped in the definition of the tree. This Bernstein-Orlicz
norm then turns up in the constants 
(\ref{gamma1*.equation}) and (\ref{gamma2*.equation}) which appear in
the generic chaining bound of Theorem \ref{generic}

\begin{definition}\label{generic-tree.definition}
Let $\delta >0$ be a constant.
A $\delta$ finite generic tree chain for $\{ Z_{\theta} \}$ is a finite labeled
tree $({\cal T} ,{\cal W} )$
such that one can apply generic chaining of $\{ Z_{\theta } \}$ along the
tree $({\cal T} , {\cal W})$ with approximation error
$\delta$. That is, 
for each $\theta \in \Theta$ there is an end  node $j_S \in
G_S$ such that the branch $\{ j_0 , \ldots , j_S \}$ satisfies
$$| Z_{\theta}| \le   \sum_{s=0}^S | W_{j_s} | + \delta . $$
\end{definition}

Let 
$({\cal T} ,{\cal W} )$ a finite labeled tree. For each end node
$k \in G_S$, we let 
$$\{ j_0 (k) , \ldots , j_S(k)\} $$ 
be the corresponding branch
(so that $j_S (k) =k $), and
we write
$$ W_s(k):= W_{j_s(k)} , \ k \in G_S , \ s=0, 1 , \ldots , S . $$
Fix a 
sequence of positive constants ${\cal L} := \{ L_s \}_{s=0}^S $.
We write for $k \in G_S$, 
\begin{equation}\label{gamma1*.equation}
\gamma_{1,*} (k) :=  \sum_{s=0}^S \| W_s (k) \|_{\Psi_{L_s}} 
\sqrt {\log (1+ | G_s | ) } , 
\end{equation}
\begin{equation}\label{gamma2*.equation}
\gamma_{2,*} (k) := \sum_{s=0}^S 
\| W_s (k) \|_{\Psi_{L_s}} { L_s } \log (1+ | G_s | ) ,
\end{equation}
$$\gamma_*(k)  := \gamma_{1,*} (k)+ { \gamma_{2,*}(k)  \over 2} . $$
Moreover, we let
$$\gamma_{1,*} := \max_{k \in G_S} \gamma_{1,*} (k) , \
\gamma_{2,*} := \max_{k \in G_S} \gamma_{2,*} (k) , \ \gamma_{*} := \max_{k \in G_S} \gamma_{*} (k) , $$
and
$$\tau_* := \max_{k \in G_S } \sum_{s=0}^S  \| W_s (k) \|_{\Psi_{L_s}} \sqrt {1+s}   , $$
and
$$L_* \tau_* := \max_{k \in G_S } \sum_{s=0}^S 
\| W_s (k) \|_{\Psi_{L_s}} (1+s) L_s .$$

\begin{theorem}\label{generic}
Let $({\cal T} ,{\cal W} )$ be a $\delta$ finite generic tree chain for
$\{ Z_{\theta} \} $. Then 
$$\PP \biggl ( \sup_{\theta \in \Theta} | Z_{\theta} | \ge \gamma_* + \delta + \tau_* \biggl [1+ {L_* \over 2}
\biggr ]  +  \tau_* \biggl [ \sqrt {t} + {L_* t \over 2}   \biggr ]
\biggr ) \le 2 \exp[-t] . $$
\end{theorem} 

\begin{remark}\label{minimize-tree2.remark}
The result of Theorem \ref{generic} may again be optimized over  all 
finite generic trees. 
\end{remark}

{\bf Proof of Theorem \ref{generic}.}
Define for $s=0 , \ldots , S$, 
$$\alpha_s  :=   \biggl [   \sqrt {\log (1+ |G_s| ) } +
  {L_s \over 2}  \log (1+ |G_s| ) \biggr ] 
  +   \biggl [ \sqrt {(1+s)(1+t) } + {(1+s)(1+t)  L_s  \over 2}\biggr ] . $$
   Using Lemma \ref{deviation-Prob.lemma}, we see that
   \begin{equation}\label{single-s.equation}
  \PP \biggl ( \max_{j \in G_s } {  |W_j |  \over
   \| W_j \|_{\Psi_{L_s} } } \ge \alpha_s  \biggr ) \le 2  \exp[-(1+t)(1+s) ] , \ 
   s=0 , \ldots , S . 
   \end{equation}
 We have
$$\PP \biggl (\max_{k } \sum_{s=0}^S | W_{s} (k)|  \ge \max_{k }
\sum_{s=0}^S\| W_s(k) \|_{\Psi_{L_s}} \alpha_s  \biggr ) $$
$$ \le \PP \biggl (\exists k : \ \sum_{s=0}^S | W_{s} (k)|  \ge \sum_{s=0}^S 
\| W_s(k) \|_{\Psi_{L_s}} \alpha_s  \biggr )  $$
$$ \le \sum_{s=0}^{S} \PP \biggl (\exists k : \   |W_{s} (k)|  \ge  \| W_s(k) \|_{\Psi_{L_s}}
\alpha_s
\biggr )   $$
$$ = \sum_{s=0}^S \PP\biggl  (\max_k 
{ | W_{s} (k)| \over \| W_s(k) \|_{\Psi_{L_s}}  }\ge 
\alpha_s   \biggr ) \le 
\sum_{s=0}^S \PP \biggl ( \max_{j \in G_s } { |W_j | \over
\| W_j \|_{\Psi_{L_s}} } \ge \alpha_s  \biggr ) . $$
Now insert (\ref{single-s.equation}) to find
$$\PP \biggl (\max_k  \sum_{s=0}^S |W_{s} (k) | \ge \max_k 
\sum_{s=0}^S  \| W_s (k) \|_{\Psi_{L_s}}  \alpha_s  \biggr ) \le
2 \sum_{s=0}^S \exp[-(1+t)(1+s) ]$$
$$ \le {2 {\rm e} ^{-(1+t)} \over
1- {\rm e} ^{-(1+t) } } \le {2 {\rm e} ^{-1} \over
1- {\rm e} ^{-1 } } \exp[-t] \le 2 \exp[-t] . $$
  We have by definition 
  $$ \max_k \sum_{s=0}^S   \| W_s(k) \|_{\Psi_{L_s}}  \biggl [ \sqrt {\log (1+ |G_s| ) } + {L_s \over 2} \log (1+ |G_s| )   \biggr ]=\gamma_*, $$
$$\max_k \sum_{s=0}^S   \| W_s(k) \|_{\Psi_{L_s}} \sqrt {(1+s)}  =\tau_*  , $$
and
$$\max_k  \sum_{s=0}^S \| W_s(k) \|_{\Psi_{L_s}} (1+s) L_s  =\tau_* L_*.$$
Therefore,
$$\max_k \sum_{s=0}^S  \| W_s(k) \|_{\Psi_{L_s}}\alpha_s \le   
\gamma_* + \tau_*  \sqrt {1+t} +   { \tau_* (1+t) L \over 2}  $$ $$
\le \gamma_* +  \tau_* + {\tau_* L_* \over 2}   + \tau_* \biggl [ \sqrt {t} + {L_*t \over 2}  \biggr ]  .$$
\hfill $\sqcup \mkern -12mu \sqcap$

Note that the constants $L_*$ and
$\tau_*$ possibly depend on the complexity of $\Theta$ through the quantities
$\{ \| W_s (k) \|_{\Psi_{L_s} } : \ k \in G_s , \ s=0, \ldots , S \} $.
Moreover, the choice of the constants ${\cal L} = \{ L_s \}_{s=0}^S $
may also depend on the complexity of $\Theta$. 
In the application to the
empirical process (see Section \ref{adaptivetruncationsection}), the latter will be indeed the case. 
We will nevertheless derive there a deviation inequality
where we put the dependency on the complexity of
$\Theta$ in the shift.

As a simple corollary of Theorem \ref{generic}, one obtains a deviation inequality in the Bernstein-Orlicz norm.
We state this for completeness. 
In Section \ref{adaptivetruncationsection} we will not apply Corollary \ref{deviation-Psi-sup.corollary} directly, because as such, it does not allow us to put
all dependency on the complexity of $\Theta$ in the shift.

\begin{corollary}\label{deviation-Psi-sup.corollary}
Let the conditions of Theorem \ref{generic} be met. Then
the combination of  this theorem with Lemma 
\ref{Prob-Psi.lemma} gives
$$\biggl \| \biggl ( \sup_{\theta \in \Theta} |  Z_{\theta} | - ( \gamma_*  + \delta +
\tau_* [1+ L_* / 2 ] )  \biggr )_+   \biggr \|_{\psi_{\sqrt 3 L_*}} \le  \sqrt 3 \tau_*  . $$
By Jensen's inequality, we then get
$$\EE \sup_{\theta \in \Theta} |  Z_{\theta} |  \le  \gamma_*  + \delta +
\tau_* \biggl [1+ { L_*  \over  2} \biggr  ] +  \sqrt 3 \tau_* \biggl [  
\sqrt {\log 2} + {\sqrt 3 L_* \over 2}  \log 2 \biggr ]  . $$
\end{corollary}

\begin{example} \label{Talagrand.example}
In  \cite{talagrand2005generic}, the sizes $|G_s|$ of generation $s$
is fixed to be 
$$| G_s |= 2^{2^{2s}} , \ s=0 , \ldots , S . $$
In that case,
$$\log (1+ |G_s | ) \le (2^{2s} +1 ) \log 2 \le 2^{2s+1}  \le 2^{2(s+1)}. $$
Hence
$$\gamma_* \le 2 \gamma_0 , $$
where 
$$ \gamma_0 := \max_{k \in G_S} \gamma_0 (k) , 
$$
and for $k \in G_S$, 
$$
 \gamma_{0} (k) :=
\gamma_{1,0} (k) + { \gamma_{2,0} (k)  \over 2} , $$
and
$$\gamma_{1,0} (k) := \sum_{s=0}^S \| W_s (k) \|_{\Psi_{L_s}}   2^s ,\ 
\gamma_{2,0} (k) :=   
  \sum_{s=0}^S  \| W_s(k) \|_{\Psi_{L_s}}   L_s 2^{2s}    . $$
  Furthermore, since
$1+s \le 2^{2s} $ for all $s \ge 0$, 
$$\tau_* \le \gamma_{1,0} := \max_{k \in G_S} \gamma_{1,0} (k) ,  $$
and
$$\tau_* L_* \le \gamma_{2,0} := \max_{k \in G_S} \gamma_{2,0} (k)  . $$
Hence, 
$$  \gamma_* + \tau_* \biggl [1+ {L_* \over 2}\biggr ] \le
3 \biggl [ \gamma_{1,0} + {\gamma_{2,0} \over 2} \biggr ] ,  $$
and
$$\sqrt 3\tau_* \biggl [  \sqrt {\log 2} + \sqrt {3 L_* \over 2} \log 2 \biggr ] 
\le \sqrt {3 \log 2} \ \gamma_{1,0} + { 3 \log 2 \over 2} \gamma_{2,0} . $$
It follows from Corollary
\ref{deviation-Psi-sup.corollary} that
 $$\EE  \sup_{\theta \in \Theta} |  Z_{\theta} | \le
(3+ \sqrt {3\log 2} ) \gamma_{1,0} + {3+ 3 \log 2 \over 2}   \gamma_{2,0} .$$
Thus, we arrive at a special case of Theorem 1.2.7 in
 \cite{talagrand2005generic}. The latter book
does not treat deviation inequalities.
\end{example}

   When using a $(\delta, \tau , {\cal L} )$ finite  tree chain,
   one takes $\| W_s (k) \|_{\Psi_{L_s}}  \le \tau 2^{-s} $ for all $s$ and $k \in G_s$. 
   In that case, the constants $\tau_*$ and $L_*$ in the bounds given
   in Corollary \ref{deviation-Psi-sup.corollary} only depend
   on the scale parameter $\tau$ and on the constants ${\cal L} = \{ L_s \}_{s=0}^S $.
   This is detailed in the next theorem. 
   
\begin{theorem} \label{Sup-Probability.theorem}
Let the conditions of Theorem \ref{Tree.theorem} be met, and define
$$\gamma :=   \tau  \sum_{s=0}^S2^{-s}  \biggl [  \sqrt {\log (1+ |G_s| ) } + {L_s \over 2} \log (1+ |G_s| )   \biggr ]  , $$
and
$$ 
L:= 
\sum_{s=0}^S { 2^{-s} L_s (1+s) \over 4}. $$
Then for all $t>0$
$$\PP \biggl ( \sup_{\theta \in \Theta } |Z_{\theta} | \ge \gamma    + \delta +4\tau
\biggl [ 
1 + {L \over 2}  \biggr ]  + 4 \tau \biggl [ \sqrt {t} + {Lt \over 2}   \biggr ]\biggr ) 
\le 2 \exp[-t] . $$

\end{theorem}

{\bf Proof of Theorem \ref{Sup-Probability.theorem}.} 
This follows from Theorem \ref{generic}, where one takes
$$\| W_s (k) \|_{\Psi_{L_s}}  \le \tau 2^{-s} . $$
  We have
$$\tau_* /\tau \le  \sum_{s=0}^S 2^{-s} \sqrt {(1+s)}  =2 \sum_{s=1}^{S+1} 2^{-s} \sqrt {s} 
 \le 2 \int_0^{\infty} 2^{-x} \sqrt x dx = {\sqrt \pi \over (\log 2)^{3/2} }   \le 4 . $$
 Moreover,
 $$L_* \tau_* \le \sum_{s=0}^S 2^{-s} L_s (1+s) = 4L . $$
\hfill $\sqcup \mkern -12mu \sqcap$ 



\section{Application to empirical processes} \label{adaptivetruncationsection}
Let ${\cal X}$ be some measurable space, and 
consider independent ${\cal X}$-valued random variables $X_1 , \ldots , X_n$.
Let ${\cal G}$ be a collection of real-valued functions on ${\cal X}$.

Write
$$P_n g := {1 \over n} \sum_{i=1}^n g(X_i), \ Pg := {1 \over n} \sum_{i=1}^n \EE g(X_i), $$
and 
$$\| g \|^2 := {1 \over n} \sum_{i=1}^n \EE g^2 (X_i) . $$
We assume the normalization
$$\sup_{g \in {\cal G} } \| g \| \le 1 . $$
We study the supremum of the empirical process
$ \{ \nu _n (g) : \ g \in {\cal G} \} $, 
where
$ \nu_n (g) := \sqrt n (P_n - P)g $. 

We recall the deviation inequality of \cite{Massart:00a}, which refines the constants in \cite{talagrand1996new}.

\begin{theorem}\label{massart.theorem} (\cite{Massart:00a}) Suppose that
for a constant $K$
\begin{equation}\label{supnorm.condition}
\sup_{g \in {\cal G} } \sup_{x \in {\cal X} } | g (x) | \le K . 
\end{equation}
Then for all $\epsilon>0$ and all $t>0$, it holds that
\begin{equation}\label{massart}
\PP \left ( \sup_{g \in {\cal G}} | \nu_n (g) | \ge 
(1+ \epsilon) \EE \sup_{g \in {\cal G}} | \nu_n (g) |+ \sqrt{2 \kappa t} +
\kappa (\epsilon) K t /\sqrt {n} \right ) \le \exp[-t] , 
\end{equation}
where $\kappa $ and $\kappa (\epsilon)$ can be taken equal to
$\kappa=4$ and $\kappa(\epsilon)= 2.5+ 32/\epsilon$. 
\end{theorem}

For the i.i.d.\ case, \cite{Bousquet:02} obtained constants remarkably
close those to for the case where ${\cal G}$ is a singleton.
In fact, \cite{Massart:00a} and \cite{Bousquet:02} and others have derived
{\it concentration} inequalities which in addition to upper bounds show similar
lower bounds for  the supremum of the empirical
process. This is complemented in \cite{Lederer11} to moment concentration inequalities
assuming only moment conditions on the envelope
$\Gamma ( \cdot ) := \sup_{g \in {\cal G}} |g(\cdot ) |$, instead of the boundedness assumption
(\ref{supnorm.condition}).

In this paper, we provide a deviation inequality of the same spirit as in the above
Theorem \ref{massart.theorem}, 
where we replace condition (\ref{supnorm.condition}) by a weaker Bernstein
condition (see (\ref{Uniform-Bernstein.condition})),
which essentially requires that the $g(X_i)$ have sub-exponential
tails, and where we also
present a deviation result in Bernstein-Orlicz norm.
These deviation results in probability and
in Bernstein-Orlicz norm are given in Theorem \ref{deviation-empirical-process.theorem}. 
We have not tried to optimize the constants. Moreover,
we replace the expectation $\EE \sup_{g \in {\cal G} }| \nu_n (g) |$ 
in (\ref{massart}) by
the upper bound we obtain from chaining arguments\footnote{This
upper bound can be shown to be (up to constants) tight in certain examples.
The upper bound following from generic chaining is modulo
constants tight for the 
general sub-Gaussian case.}.
Deviation inequalities
for the sub-exponential case can be found in literature (see 
e.g.\ 
\cite{viens2007supremum}), but these do not cover the more refined
interpolation of sub-Gaussian and sub-exponential tail behavior. 
The above cited work also contains lower bounds for suprema,
thus completing the results to concentration inequalities.

Now our first aim is to show that entropy with 
bracketing conditions allow one to construct a finite tree chain.
We recall here the definition of a bracketing set and entropy with bracketing
(see \cite{Blum:55}, or see \cite {vanderVaart:96}, 
  \cite{vandeGeer:00} and their references).

\begin{definition} Let $s>0$ be arbitrary. A $2^{-s}$-bracketing set for 
$\{ {\cal G}, \| \cdot \| \} $
is a finite collection of functions $\{ [\tilde g_j^L, \tilde g_j^U ] \}_{j=1}^{\tilde N_s} $ satisfying
$\| \tilde g_j^U - \tilde g_j^L \| \le 2^{-s} $ for all $j$, and such that for
each $g \in {\cal G}$ there is a $j \in \{ 1 , \ldots , \tilde N_s \} $ such that
$\tilde g_j^L \le g \le \tilde  g_j^U $.  If no such finite collection exists, we write
$\tilde N_s = \infty$. 
\end{definition}

We also introduce a generalized bracketing set, in the spirit of 
\cite{vandeGeer:00}.

\begin{definition} Let $K>0$ be a fixed constant. A generalized bracketing set for 
$ {\cal G} $
is a finite collection of functions $\{ [\tilde g_j^L,\tilde  g_j^U ] \}_{j=1}^{\tilde N_0} $ satisfying
for all $j$
$$P| \tilde g_j^U - \tilde g_j^L |^m \le {m! \over 2} (2K)^{m-2}  , \ 
m=2,3, \ldots , $$  and such that for
each $g \in {\cal G}$ there is a $j \in \{ 1 , \ldots , { \tilde N}_0 \} $ such that
$\tilde g_j^L \le g \le \tilde g_j^U $.  Write $\tilde N_0= \infty$ if no such finite
collection exists. 
\end{definition}

A special case is where the envelope function 
$\Gamma := \sup_{g \in  {\cal G}} | g|  $
satisfies the Bernstein condition
$$P\Gamma^m \le {m! \over 2} (2K)^{m-2}  , \ 
m=2,3, \ldots .$$ Then one can take
$[ - \Gamma, \Gamma]$ as generalized bracketing set, consisting of only one
element. 

In what follows, we let for each $s \in \Nat$, 
$\tilde  N_s$ be the cardinality of a minimal $2^{-s}$-bracketing set for ${\cal G}$.
The $2^{-s} $-entropy with bracketing of ${\cal G}$ is
$$\tilde H_s := \log (1+ \tilde N_s) , \ s \in \Nat. $$ Moreover, $\tilde N_0 $ is
the cardinality of a minimal generalized bracketing set, and we let
$${ \tilde H}_0 := \log (1+ { \tilde N}_0). $$Finally, we write
\begin{equation}
\label{product-generations}
N_s :=  \prod_{k=0}^s \tilde N_k , \
H_s := \log (1+ N_s ) , \ s \in \Nat_0 . 
\end{equation} 

The following theorem uses arguments of  \cite{Ossiander:87},
and is comparable to Theorem 2.7.11 in
\cite{talagrand2005generic} (who adapts the technique
of \cite{Ossiander:87}).
However, we do not use generic chaining here. On the other hand, our results
lead to the more involved 
deviation inequalities as given in Theorem \ref{deviation-empirical-process.theorem}.

\begin{theorem} \label{tree-bracketing.theorem}
Suppose that for some constant $K \ge 1$, one has the
Bernstein condition
\begin{equation}\label{Uniform-Bernstein.condition}
\sup_{g \in {\cal G}}  P|  g|^m \le {m! \over 2} K^{m-2}  , \ 
m=2,3, \ldots .
\end{equation}
Let $S$ be some integer, $\tau:= 3 \sqrt {6}$ and $\delta: =
4 \sqrt n \sum_{s=1}^{S}   2^{-{2s} } / K_{s-1} +\sqrt n 2^{-S}  $,
where $\{K_{s-1} \}_{s=1}^S $ is an arbitrary deceasing sequence of positive constants
(called truncation levels). 
 Suppose that $\tilde N_s < \infty$ for all
$s=0, \ldots , S$. 
Then there is a $(\delta , \tau, {\cal L} )$ finite tree chain for
$\{ \nu_n (g) \}$, with $|G_s | \le N_s$, $s=0, \ldots , S$, and with
$$L_0 =  { 4 \sqrt 6   K \over \sqrt n} , \ 
L_s =  { 2  \sqrt 6 \ 2^s  K_{s-1} \over  3 \sqrt {n} }  , \ s=1 , \ldots , S .  $$
\end{theorem} 

As a consequence, we can derive a bound for the expectation
of the supremum of the empirical process.

\begin{theorem} \label{expectation-empirical-process.theorem}
Assume the Bernstein condition (\ref {Uniform-Bernstein.condition}).
Let   
$$\bar {\bf E}_S := 
2^{-S} \sqrt {n} + 14 \sum_{s=0}^S 2^{-s} \sqrt {6 \tilde H_s } 
+ 6^2 K {  \tilde H_0   \over \sqrt {n} } . $$
Then one has
$$\EE \biggl ( \sup_{g \in {\cal G} } |\nu_n (g) | \biggr ) \le \min_{S } 
\bar {\bf E}_S . $$
\end{theorem}

\begin{remark} \label{S=0.remark} When $\Theta $ is finite, say $|\Theta |=p$,
one may choose a bound with $S=\delta=0$, and
$\tilde H_0 \le \log (1+ p)$.  Theorem then
\ref{expectation-empirical-process.theorem} yields - up to constants - 
the same bound as in (\ref{S=0}). 
\end{remark}

Finally, we present the main result of this section.
We give deviation results in probability and
in Bernstein-Orlicz norm, where the dependency on the complexity
of ${\cal G}$ is only in the shift.

\begin{theorem} \label{deviation-empirical-process.theorem}
 Assume the Bernstein condition
(\ref {Uniform-Bernstein.condition}). 
Define as in Theorem \ref{expectation-empirical-process.theorem}, 
$$\bar {\bf E}_S := 
2^{-S} \sqrt {n} + 14 \sum_{s=0}^S 2^{-s} \sqrt {6 \tilde H_s } 
+ 6^2 K {  H_0   \over \sqrt {n} } .$$
Let
$$ \tilde L:= {  \sqrt 6 K \over 2 \sqrt n } . $$
Then for all $t>0$, 
$$\PP \biggl ( \sup_{g \in {\cal G}} |  \nu_n (g) | \ge \min_S \bar { \bf E}_S +
 6^2 K / \sqrt n + 24 \sqrt 6 \biggl [ \sqrt t + { \tilde L t \over 2} \biggr ]
  \biggr ]  \biggr ) \le 2 \exp[-t] .  $$
  Moreover,
  $$\biggl \| \biggr (  \biggl [ \sup_{g \in {\cal G}} |  \nu_n (g)| \biggr ] - \biggl [ \min_{S} \bar {\bf E}_S  + 6^2 K /\sqrt n + 24 \sqrt 6 \biggr ] 
 \biggr)_+  \biggr  \|_{\Psi_{\sqrt 3 \tilde L}  } \le 72 \sqrt 2   .$$
\end{theorem}

Theorem \ref{deviation-empirical-process.theorem}
can be compared to results in \cite{Adamczak:08}. One sees that
our bound replaces the sub-exponential Orlicz-norm
$$\biggl \| \max_{1 \le i \le n} \sup_{g \in {\cal G} } | g(X_i) |\biggr  \|_{\Psi} , \
\Psi (z) = \exp (z) -1, z \ge 0 , $$
occurring in \cite{Adamczak:08} by a constant proportional to $K$,
which means we generally gain a $\log n $-term.
On the other hand, the shift in \cite{Adamczak:08} is
up to a factor $(1+\epsilon)$ equal to the expectation 
$$\EE \sup_{g \in {\cal G} } | \nu_n (g) | , $$
as in  \cite{Massart:00a}) (whose result is cited here
in Theorem \ref{massart.theorem}). 

\begin{remark} Again, when $|\Theta |= p$ is finite, one can choose $S=\delta=0$,
and $\tilde H_0 \le \log (1+p)$. 
as in Remark \ref{S=0.remark}. Theorem
\ref{deviation-empirical-process.theorem} then reduces to the usual  union bound type deviation
inequalities for the maximum of finitely many random variables (that is,
the results are - up to constants -
a special case of Lemmas \ref{deviation-Prob.lemma}
and \ref{deviation-Psi.lemma}).
\end{remark}

 \section{Proofs for Section \ref{adaptivetruncationsection}}
 \label{proofs.section} 
 
\subsection{Proof of Theorem \ref{tree-bracketing.theorem}}
This follows from similar arguments as in \cite{vandeGeer:00},
who uses in turn ideas of \cite{Ossiander:87}. 
Let for $s=1,  \ldots , S$, 
$$\{ [ \tilde g_j^{s,L} , \tilde g_j^{s,U} ]\}_{j=1}^{\tilde N_s} $$
be a minimal $2^{-s }$-bracketing set for $\| \cdot \| $.
Let
$\{ [\tilde g_j^{0,L}, \tilde g_j^{0,U} ] \}_{j=1}^{{\tilde N}_0} $ be a generalized bracketing set.

Consider some  $g\in {\cal G}$, and  let $
[\tilde g^{0,L} , \tilde g^{0,U} ] $ be the corresponding
generalized bracket, and for all $s\in \{1 , \ldots , S \} $, let the corresponding brackets
be
$[\tilde g^{s,L} , \tilde g^{s,U} ]$. Thus
$$\tilde g^{s,L} \le g \le \tilde g^{s,U} , \ s=0, \ldots , S, $$
and
$$P | \tilde g_{0,U} - \tilde g_{0,L} |^m \le { m! \over 2} (2K)^{m-2}, \ m=2,3, \ldots, $$ $$
P | \tilde g^{s,U} - \tilde g_{s.L} |^2\le 2^{-2s} , \ s=1, \ldots , S . $$
If for some $s$ there are several brackets in $\{ [ \tilde g_j^{s,L}, \tilde g_j^{s,U} ] \}_{j=1}^{\tilde N_s}$
corresponding to $g$, 
we choose a fixed but otherwise arbitrary one.
Define
$$g^{s,L} := \max_{ 0 \le k \le s} \tilde g^{k,L} , \ g^{s,U} := \min_{ 0 \le k \le s} \tilde g^{k,U}.$$
Then 
$$g^{0,L} \le g^{1,L} \le \cdots  \le g^{S,L} \le g \le g^{S,U} \le
\cdots \le g^{1,U} \le g^{0,U} , $$ and 
moreover $ g^{s,U} - g^{s,L} \le \tilde g^{s,U} - \tilde g^{s,L} $.
Denote the difference between upper and lower bracket by
$$\Delta^s := g^{s,U}-g^{s,L} , \ s=0 , \ldots , S.$$
The differences
$\Delta^s$ are decreasing in $s$.
Furthermore, $\| \Delta^s \|  \le 2^{-s} $, for all $s\in \{ 0,1 , \ldots , S \} $.

Let ${\cal N}_s := | \{ [ g_j^{s,L}, g_j^{s,U}] \}  | $, $s=0, \ldots , S$. 
It is easy to see
that
$${\cal N}_s \le \prod_{k=0}^s \tilde N_k =: N_s, \ s=0, \ldots , S. $$
We define a tree with end nodes
$\{ 1 , \ldots , {\cal N}_S\}$. At each end node $j$ sits a pair
of brackets $[g_j^{S,L} , g_j^{S,U} ]$.
For each $s=0 , \ldots , S-1$, we define the parents at
generation $s$ as follows.
Let
$$\tilde V_k^s := \{ l: [g_k^{s-1,L}, g_k^{s-1,U} ] \ {\rm forms \ a} \ 2^{-(s-1)}{\rm -bracket \  for} \ 
[ g_l^{s,L} , g_l^{s,U} ] \}. $$
Then $\cup_{k=1}^{{\cal N}_{s-1} }  \tilde V_k^s = \{ 1 , \ldots , {\cal N}_s\}$, that is, for  each bracket
$[ g_l^{s,L}, g_l^{s,U} ] $ there is a $k\in \{ 1 , \ldots , {\cal N}_{s-1}\}$ with $
l \in \tilde V_k$.  To see this, we note that for each $l$, there is a function
$g$ with $g_l^{s,L} \le g \le g_l^{s,U}$, and  by the above construction, there
is a $k$ with $g_k^{{s-1}, L} \le g_l^{s, L} \le g \le g_l^{s,U} \le g_k^{s-1, U} $.
We let $\{ V_k^s \}_{k=1}^{{\cal N}_{s-1}} $ be a disjoint version
of $\{ \tilde V_k^s \} $, e.g., the one given by
$$V_1^s = \tilde V_1^s , \ V_k^s = \tilde V_k^s \backslash \cup_{l=1}^{k-1} \tilde V_l^s, \ k=1 , \ldots , {\cal N}_{s-1} . $$
We let
$${\rm parent} (j_s) = k \ {\rm  if} \  j_s \in V_k^s.   $$

We now turn to an adaptive truncation device.
For for each $s=0, \ldots , S-1$, we are given 
truncation levels $K_s$, such that
$K_s$ is assumed to be decreasing in $s$. 
Let $g$ be fixed and
$$g^{0,L} \le g^{1,L} \le \cdots  \le g^{S,L} \le g \le g^{S,U} \le
\cdots \le g^{1,U} \le g^{0,U} . $$
Define 
$$\Delta^s := g^{s,U} - g^{s,L}, \ y_s := {\rm l} \{ \Delta^s \ge K_s \} . $$

Then
$$K_s {\rm l} \{ y_s=1\} \le \Delta^s {\rm l}  \{ y_s=1  \}
 , \ s=0 , \ldots , S-1 , $$
which implies (for $s=0, \ldots , S-1$)
$$P\Delta^s {\rm l}  \{ y_s=1 \}  \le
{ P| \Delta^s|^2 \over K_s } \le { 2^{-2s}   \over K_s} $$

We can write any $g \in {\cal G}$ as
\begin{eqnarray}\label{telescope.equation}
 & &g=
\sum_{s=1}^S ( g - g^{0,s} ) {\rm l} \{ y_s =1 , y_{s-1} = \ldots = y_0 =0 \} \\ &&+
\sum_{s=1}^S ( g^{s,L} - g^{s-1, L} ) {\rm l} \{ y_{s-1}= \ldots = y_0 =0 \} +
g_{0,L} + ( g- g^{0,L} ) {\rm l} \{ y_0=1 \}
\nonumber
\end{eqnarray} 
Let
$$W_{j_0 } := |\nu_n ( g^{0,L})  |+ |\nu_n (\Delta^0 ) | , $$
$$ W_{j_s }:=  |\nu_n (\Delta^s {\rm l} \{ y_{s-1}=0 \}
) |+
|\nu_n ((g^{s,L} - g^{s-1, L} ){\rm l}\{ y_{s-1}=0   \} ) | , \ 
s=1 , \ldots , S . $$
Then it follows from (\ref{telescope.equation}) that
$$| \nu_n (g) | \le    \sum_{s=0}^S | W_{j_s} | + \sqrt {n}
\sum_{s=0}^S P \Delta^s {\rm l} \{ y_s=1 \}
 \le   \sum_{s=0}^S | W_{j_s} | + \delta, $$
for 
$$\delta =
 \sqrt {n} \sum_{s=1}^{S} { 4 \ 2^{-2s} \over K_{s-1} }+  \sqrt {n} 2^{-S}  . $$

Note now that
$$( P |g^{0,L} |^m)^{1/m}  \le 
( P |g |^m)^{1/m} + ( P |\Delta^0 |^m)^{1/m} $$ $$ \le  ( { m! \over 2} K^{m-2} )^{1/m}
+  ( { m! \over 2} (2K)^{m-2} )^{1/m} \le 2 ( { m! \over 2} (2K)^{m-2} )^{1/m} , $$
so
$$ P |g^{0,L} |^m \le {m!\over 2} (4K)^{m-2} 2^2 . $$
By Corollary \ref{Bernstein-Psi.corollary}
$$ \| \nu_n ( g^{0,L} ) \|_{\psi_{L_0} } \le 2 \sqrt 6 ,$$
for
$$L_0 = {  \sqrt 6  }  (8K /2)/\sqrt n  = 4 \sqrt 6 {/ \sqrt n} , $$
where we  multiplied by a factor 2 because the Bernstein
condition for the centered functions holds with the above $4K$ replaced by $8K$.
Moreover, 
$L_0 = {\sqrt 6 } (4K)/\sqrt n$, so again by Corollary \ref{Bernstein-Psi.corollary},
$$\| \nu_n ( \Delta^0 ) \|_{\Psi_{L_0} } \le \sqrt 6 . $$
The triangle inequality gives 
$$\biggl \| | \nu_n (g^{0,L} ) | + |\nu_n (\Delta^0) |  \biggr \|_{\Psi_{L_0} } \le 3 \sqrt 6 =:
\tau  . $$

Moreover, for $s=1 , \ldots , S$, 
$$| ( g^{s, L} - g^{s-1, L}) {\rm l} \{  y_{s-1}=0 \}   | \le \Delta^{s-1} \le
K_{s-1}    , \  \| \Delta^{s-1} \| \le  2^{-s+1} , $$
and 
$$  \Delta^s  
{\rm l}\{   y_{s-1}=0 \} \le  \Delta^{s-1} \le K_{s-1}  , \ \| \Delta^s \| \le 2^{-s}. $$
So, again by  Corollary \ref{Bernstein-Psi.corollary}, we may take 
$$L_s :=   \sqrt {6}  \  2^s \max( {2 \over 3} K_{s-1} /2, {2 \over 3} K_{s-1} )/\sqrt n =
{2  \sqrt {6}K_{s-1} \over 3 \sqrt n }  , 
\ s=1 , \ldots , S .  $$
Then, again by the triangle inequality, 
$$\biggl \| | \nu_n (( g^{s, L} - g^{s-1, L}) {\rm l} \{ y_{s-1}=0 \} ) |+
|\nu_n (\Delta^s  {\rm l} \{  y_{s-1}=0 \}) | \biggr \|_{\Psi_{L_s}} \le 3 \sqrt 6 \ 2^{-s} . $$

\hfill $\sqcup \mkern -12mu \sqcap$

\subsection{Three technical lemmas}

To apply the result of Theorem \ref{tree-bracketing.theorem}, we need three technical lemmas.
First we need a bound for $N_s:= \prod_{k=0}^s \tilde N_s$,
or actually for $H_s := \log (1+ N_s)$.

\begin{lemma} \label{5decreasingcoveringslemma}
Let $s \in \{ 0, \ldots , S\}$, $H_s := \log (1+ \prod_{k=0}^s \tilde N_k ) $ and
$ \tilde H_s := \log (1+ \tilde N_s)$.
It holds that
$$\sum_{s=1}^S 2^{-s} \sqrt {H_s} \le \sqrt {\tilde H_0} + 2
\sum_{s=1}^S 2^{-s} \sqrt {\tilde H_s} . $$
\end{lemma}

{\bf Proof of Lemma \ref{5decreasingcoveringslemma}.} 
We have
$$\sqrt {H_s} \le \sum_{k=0}^s \sqrt {\tilde H_k } , $$
so
$$\sum_{s=1}^S 2^{-s} \sqrt {H_s} \le
\sum_{s=1}^S2^{-s}  \sqrt {  \tilde H_0 } +\sum_{s=1}^S 2^{-s} \sum_{k=1}^s \sqrt 
{\tilde H_k } $$
$$\le  \sqrt {\tilde H_0} + \sum_{k=1}^S \sum_{s=k}^S 2^{-s}  \sqrt {\tilde H_k } 
 \le  \sqrt {\tilde H_0} +2  \sum_{k=1}^S 2^{-k}   \sqrt {\tilde H_k } .$$
\hfill $\sqcup \mkern -12mu \sqcap$

The next lemma inserts a special choice for the truncation levels
$\{ K_s \}$, and then establishes a bound for the expectation of the supremum
of the empirical process, derived from the one of Theorem 
\ref{Tree.theorem}.

\begin{lemma}\label{1ChoiceK.lemma} Let
Let $S$ be some integer and $\epsilon \ge 0$ be an arbitrary constant.
Take
$$K_{s-1} := 2^{-s}   \sqrt n \biggl (  {  \sqrt 6   \over 3 \sqrt { \log (1+ N_s) }}
\wedge{ 1 \over \epsilon}  \biggr )   , \ s=1 , \ldots , S,$$
where $u\wedge v$ denotes the minimum of $u$ and $v$. 
Define as in Theorem \ref{tree-bracketing.theorem}, 
$$L_0 :=  {4 \sqrt 6   K \over \sqrt n} , \ L_s :=  { 2 \sqrt 6 \  2^s K_{s-1} \over 3 \sqrt {n} }  , \ s=1 , \ldots , S ,  $$
$$\delta: =
4 \sqrt n \sum_{s=1}^{S}   2^{-{2s} } / K_{s-1} +\sqrt n 2^{-S}  , $$
and $\tau:= 3\sqrt 6$.
Let
$${\bf E}_S := \tau  \sum_{s=0}^S 2^{-s} \biggl [ \sqrt {\log (1+ N_s)} + {L_s \over 2} 
\log (1+ N_s ) \biggr ] +  \delta  . $$
Then
$${\bf E}_S   \le \bar {\bf E}_S +4\epsilon, $$
where
$$\bar {\bf E}_S := 
2^{-S} \sqrt {n} + 14 \sum_{s=0}^S 2^{-s} \sqrt {6 \tilde H_s } 
+ 6^2 K {  H_0   \over \sqrt {n} }  . $$
\end{lemma}

{\bf Proof of Lemma \ref{1ChoiceK.lemma}.}
We have
$${\bf E}_S = \sum_{s=1}^{S} 
{4 \ 2^{-2s} \sqrt {n} \over K_{s-1}  } + 2^{-S} \sqrt {n}   +
\tau \sqrt {\log (1+ N_0)}  + 2 \sqrt 6\  \tau  K {  \log (1+ N_0)  \over \sqrt {n} } $$ $$
+ 
\tau \sum_{s=1}^S 2^{-s} \biggl [\sqrt {\log (1+ N_s)}  + 
{1\over 3} \sqrt 6  \  2^s K_{s-1}  { \log (1+ N_s ) \over \sqrt {n}}  \biggr ] $$
$$=
\sum_{s=1}^{S} 
{4\  2^{-2s} \sqrt {n} \over K_{s-1}  } + 2^{-S} \sqrt {n}   +
3  \sqrt {6 \log (1+ N_0)}  +  6^2 K {  \log (1+ N_0)  \over \sqrt {n} } $$ $$
+ 
3 \sum_{s=1}^S 2^{-s}\sqrt {6 \log (1+ N_s)}  + \sum_{s=1}^S 
6   K_{s-1}  { \log (1+ N_s ) \over \sqrt {n}} =I+II+III, $$
where 
$$I:= 2^{-S} \sqrt n + 3  \sqrt {6 \log (1+ N_0)}  +  6^2 K {  \log (1+ N_0)  \over \sqrt {n} }, $$
$$II:= 3 \sum_{s=1}^S 2^{-s}\sqrt {6 \log (1+ N_s)}, $$
and
$$III:= \sum_{s=1}^{S} 
{4 \ 2^{-2s} \sqrt {n} \over K_{s-1}  } +  \sum_{s=1}^S 
6   K_{s-1}  { \log (1+ N_s ) \over \sqrt {n}}.$$

Insert
$$K_{s-1}={1 \over 3} \sqrt 6  \ 2^{-s} \sqrt  { n  \over \log (1+ N_s) }
\wedge 2^{-s} { \sqrt n \over \epsilon}  , \ s=1 , \ldots , S.$$
Note that $K_s$ is decreasing in $s$. 
Moreover
$$ {4 \ 2^{-2s} \sqrt {n} \over K_{s-1}  } +
6  K_{s-1}  { \log (1+ N_s ) \over \sqrt {n}}\le 
4  \sqrt {6}\ 2^{-s}  \sqrt {  \log (1+ N_s ) } + 4\  2^{-s} \sqrt n \epsilon . $$

We find
$$III \le 
4  \sqrt {6}\sum_{s=1}^S 2^{-s}  \sqrt {  \log (1+ N_s ) } + 4\epsilon,$$
so that
$$II+III\le 7  \sqrt {6}\sum_{s=1}^S 2^{-s}  \sqrt {  \log (1+ N_s ) } + 4\epsilon.$$

Now apply Lemma \ref{5decreasingcoveringslemma}.
This gives
$$ II + III \le 7 \sqrt 6 \sqrt { \log (1+ \tilde N_0) }+14 \sqrt 6 
\sum_{s=1}^S 2^{-s}  \sqrt {  \log (1+ \tilde N_s ) } + 4\epsilon.$$
Hence,
$$ I+II+III\le 2^{-S} \sqrt n  + 6^2K {  \log (1+ \tilde N_0)  \over \sqrt {n} }  +10 \sqrt 6  
\sqrt {\log (1+ \tilde N_0) } $$ $$+  14 \sqrt 6 
\sum_{s=1}^S 2^{-s} \sqrt { \log (1+ \tilde N_s)}  + 4\epsilon$$
$$ \le 2^{-S} \sqrt {n} + 14  \sqrt 6 \sum_{s=0}^S 2^{-s} \sqrt {  \log (1+ \tilde N_s)}
+ 6^2K {  \log (1+ \tilde N_0)  \over \sqrt {n} } +4\epsilon. $$

\hfill $ \sqcup \mkern -12mu \sqcap$

We now derive some bounds which will be used for obtaining the deviation
inequalities in probability and in Bernstein-Orlicz norm
of Theorem \ref{deviation-empirical-process.theorem}. 

\begin{lemma}\label{2ChoiceK.lemma}
Let the constants $\{K_{s-1} \}_{s=1}^S$, $\{ L_s\}_{s=0}^S$, and $\tau$ be as
in Lemma \ref{2ChoiceK.lemma}. 
Let
$$L:= 
\sum_{s=0}^S 2^{-s} { L_s (1+s) \over 4}.  $$
Then 
$$L \le \sqrt 6 K/\sqrt n  +2 \wedge {\sqrt 6 \over \epsilon} , $$
and
$$4 \tau(1+L/2)   \le 6^2 K/\sqrt n + 24 \sqrt 6    .$$
\end{lemma}

{\bf Proof of Lemma \ref {2ChoiceK.lemma} .} 
We have
$$L= {L_0 \over 4} + \sum_{s=1}^S {2^{-s}  L_s (1+s)  \over 4 } $$ $$ =
{  \sqrt 6 K \over \sqrt n}    + \sum_{s=1}^S  {  (1+s) K_{s-1}  \over \sqrt {6n}} .$$
 But
$$\sum_{s=1}^S  2^{-s} (1+s) \le 2 \int_0^{\infty} 2^{-x} x dx = {2 \over (\log 2)^2 } , $$
and since $H_s =\log (1+ N_s ) \ge \log(2)$, 
$$K_{s-1} \le 2^{-s} \sqrt n \biggl ( {\sqrt 6 \over 3 ( \log (2))^{1/2} } \wedge {1 \over \epsilon} 
\biggr ) . $$
Hence,
$$L \le { \sqrt 6 K \over \sqrt n} + { 2 \over \sqrt 6 (\log 2)^2} 
\biggl ( {\sqrt 6 \over 3 ( \log (2))^{1/2} } \wedge {1 \over \epsilon} 
\biggr ) $$ $$= { \sqrt 6 K \over \sqrt n} + { 2 \over 3 (\log 2)^{(5/2)}} \wedge 
 { 2 \over  6 (\log 2)^2}
{\sqrt 6 \over  \epsilon}
$$ $$
   \le {\sqrt 6 K\over \sqrt n}  + 2 \wedge {\sqrt 6 \over \epsilon} .$$
As $\tau=3 \sqrt 6$, we get
$$4 \tau(1+L/2) \le 6^2 K/\sqrt n + 24 \sqrt 6  . $$
\hfill $\sqcup \mkern -12mu \sqcap$

\section{Proof of Theorems \ref{expectation-empirical-process.theorem}
and \ref{deviation-empirical-process.theorem}} \label{lastsection}

{\bf Proof of Theorem \ref{expectation-empirical-process.theorem}.}
This follows from Theorem \ref{Tree.theorem},
Theorem \ref{tree-bracketing.theorem}, and Lemma \ref{1ChoiceK.lemma}
with $\epsilon=0$.
\hfill $\sqcup \mkern -12mu \sqcap$

{\bf Proof of Theorem \ref{deviation-empirical-process.theorem}.} 
Let $t>0$ be arbitrary. 
Note that $\bar {\bf E}_S$ is as in Lemma \ref{1ChoiceK.lemma}.
Apply the bounds of Lemma  \ref{2ChoiceK.lemma} with $\epsilon=3 \sqrt {t}$ for
the constant $L$ defined there.
Then
$$\tau(4+2L) +4\epsilon + 4\tau \biggl [ \sqrt t+ {L t \over 2} \biggr ] \le
6^2 K/\sqrt n + 24 \sqrt 6 + 4 \epsilon +  12 \sqrt {6t} + 2\tau {  \sqrt 6  K t  \over \sqrt n}
+ 2\tau  {\sqrt 6 t  \over \epsilon}  $$
$$=  6^2 K /\sqrt n + 6^2 Kt/\sqrt n + 24 \sqrt 6+ 12 \sqrt{6t} +24 \sqrt t  $$
$$\le  6^2 K / \sqrt n + 24 \sqrt 6+ 24 \sqrt 6 \biggl [ \sqrt t + { \tilde L t \over 2} \biggr ] , $$
where
$$ \tilde L:= {  \sqrt 6 K \over 2 \sqrt n } . $$
Then by Theorem \ref{Sup-Probability.theorem}, 
$$\PP \biggl ( \sup_{g \in {\cal G}} |  \nu_n (g) | \ge \min_S \bar {\bf E}_S +
 6^2 K / \sqrt n +24 \sqrt 6+ 24 \sqrt 6 \biggl [ \sqrt t + { \tilde L t \over 2} \biggr ]
  \biggr ]  \biggr ) \le 2 \exp[-t] .  $$
and by Lemma \ref{Prob-Psi.lemma}
$$\biggl \| \biggr (  \biggl [ \sup_{g \in {\cal G}} |  \nu_n (g)| \biggr ] - \biggl [ \min_{S} \bar
{\bf E}_S  + 6^2 K /\sqrt n + 24 \sqrt 6 \biggr ] 
 \biggr)_+  \biggr  \|_{\Psi_{\sqrt 3 \tilde L}  } \le 72 \sqrt 2   .$$

\hfill $\sqcup \mkern -12mu \sqcap$

\bibliographystyle{plainnat}

\bibliography{reference}




\end{document}